\documentclass[10 pt]{amsart}

\usepackage{amssymb,amsmath,amsfonts,mathrsfs}

\newtheorem{theorem}{Theorem}[section]
\newtheorem{thm}{Theorem}[section]
\newtheorem{lemma}[theorem]{Lemma}
\newtheorem{lem}[theorem]{Lemma}
\newtheorem{cor}[theorem]{Corollary}

\theoremstyle{definition}

\newtheorem{example}[theorem]{Example}

\theoremstyle{remark}
\newtheorem{remark}[theorem]{Remark}
\newtheorem{rem}[theorem]{Remark}
\numberwithin{equation}{section}

\newcommand{\xcomp}{C_{c}(X)}
\newcommand{\lw}{{\ell^{2}_{\mu}}(X)}

\newcommand{\IM}{\operatorname{Im}}
\newcommand{\RE}{\operatorname{Re}}

\title[Schr\"odinger operators on graphs]{Essential self-adjointness of non-semibounded Schr\"odinger operators on infinite graphs}
\author{Ognjen Milatovic}
\address{Department of Mathematics
and Statistics \\ University of North Florida \\ Jacksonville, FL
32224 \\ USA.}
\email{omilatov@unf.edu}

\begin{document}
\maketitle

\begin{abstract} We work in the setting of infinite, not necessarily locally finite, weighted graphs. We give a sufficient condition for the essential self-adjointness of (discrete) Schr\"odinger operators $\mathcal{L}_{V}$ that are not necessarily lower semi-bounded. As a corollary of the main result, we show that  $\mathcal{L}_{V}$ is essentially self-adjoint if the potential $V$ satisfies $V(x)\geq -b_1-b_2[\rho(0,x)]^2$, for all vertices $x$, where $o$ is a fixed vertex, $b_1$ and $b_2$ are non-negative constants, and $\rho$ is an intrinsic metric of finite jump size, such that the restriction of the weighted vertex degree to every ball corresponding to $\rho$ is bounded (not necessarily uniformly bounded).
\end{abstract}

\section{Introduction}\label{S:main} Inspired by the seminal papers~\cite{Keller-Lenz-10, Keller-Lenz-09,W-08}, the analysis on (discrete) weighted graphs has experienced a significant development in the last fifteen years. A fundamental aspect of this endeavor is the study of existence and uniqueness of  self-adjoint extensions of discrete Schr\"odinger operators, a topic that has sparked the interest of several researchers, leading to quite a few papers, such as~\cite{vtt-11-2,vtt-11-3,Gol,GKS-15,HKMW,Jor-08,
Keller-Lenz-10,Keller-Lenz-09,KN-22,Milatovic-Truc,MS-18,Torki-10}. With regard to the self-adjointness of (primarily lower semi-bounded) magnetic Schr\"odinger operators (acting on vector bundles) over infinite (not necessarily locally finite) graphs, the most general results to date are contained in~\cite{MS-18} and, for (lower semi-bounded) Schr\"odinger operators acting on scalar functions on discrete (as well as metric) graphs, in the paper~\cite{KN-22}. For additional pointers to the literature on self-adjointness of discrete Laplace/Schr\"odinger operators on graphs, we refer the reader to the papers~\cite{Keller-15,KN-22,KN-23,MS-18} and the monographs~\cite{KLW-book,KN-22-book}.

On a related research track, motivated by~\cite{Masamune-09}, some authors have studied Laplace-type operators acting on 1-forms, 2-simplicial complexes, and triangulations; see, for instance,~\cite{AACTH-23,AT-15,BGJ-15,CY-17}. Recently, this line of investigation has been advanced in the paper~\cite{BK-25} concerning Hodge Laplacians on general simplicial complexes. On yet another track, in~\cite{IKMW-24}, the authors investigated the self-adjointness (and other properties) of the weighted Laplacian on birth-death chains. Interesting connections between essential self-adjointness and other properties (such as $L^2$-Liouville property) of the Laplacian on graphs (and manifolds) were explored in~\cite{HMW-21}. Additionally, far-reaching studies of (quadratic) form extensions (with applications to graphs) were conducted by the authors of~\cite{LSW-17,LSW-16}. Lastly, the authors of~\cite{ABTH-19} and~\cite{ABTH-20}, studied various properties of non-self-adjoint operators.

While there is extensive literature on the self-adjointness of lower semi-bounded discrete Schr\"odinger operators, we are familiar with just two types of results without the lower semi-boundedness assumption. As a brief reminder, in the context of locally finite graphs, the paper~\cite{Milatovic-Truc} contains a Sears-type self-adjointness result for Schr\"odinger operators with a potential $V(x)\geq -q(x)$, where $q\geq 1$ is a minorant function satisfying certain conditions, such as the completeness of a suitable path metric, obtained by modifying the usual degree-path metric  with the help of the function $q$ (see section~\ref{S:example} below for the notion of degree-path metric). Furthermore, in the setting of general (not necessarily locally finite) weighted graphs, an important self-adjointness result (whose formulation we recalled in the text of example~\ref{E:ex-3} below) was established in proposition 2.2 of~\cite{Gol}.

In the main result of the present article, theorem~\ref{T:main-1}, we give another sufficient condition for the essential self-adjointness of the Schr\"odinger operator $\mathcal{L}_{V}$,  with a potential function $V$ (see section~\ref{SS:expressions} for the description of $\mathcal{L}_{V}$). Our theorem~\ref{T:main-1} is situated on general (not necessarily locally finite) weighted graphs (with vertex set $X$), satisfying the condition $[\mathcal{L}_{V}\xcomp]\subseteq\lw$, where $\xcomp$ is the space of finitely supported functions on $X$ and $\lw$ is the space of square integrable functions with weight $\mu$, as in~(\ref{E:l-p-def}) below. We describe the remaining assumptions of theorem~\ref{T:main-1}, guaranteeing the essential self-adjointness of $\mathcal{L}_{V}|_{\xcomp}$: (i) we assume that our graph (with vertex set $X$) has an intrinsic metric $\rho$, such that the weighted vertex degree is bounded on balls $B_{\rho}(x,r)$ (centered at $x\in X$ with radius $r\geq 0$); (ii) we assume that our potential $V$ satisfies $V=U-W$, where $U\geq 0$ and $W\geq 0$, and
\[
|\nabla W|^2(x)\leq c_1+c_2|W(x)|,
\]
for all $x\in X$, where $c_1\geq 0$ and $c_2\geq 0$ are constants, and $|\nabla W|^2$ is as in~(\ref{E:grad-sq}) below.

As a consequence of theorem~\ref{T:main-1}, we get corollary~\ref{C:main-1}, establishing the essential self-adjointness of $\mathcal{L}_{V}|_{\xcomp}$ under the same hypotheses on the graph as in theorem~\ref{T:main-1} (additionally, with ``finite jump size" requirement on the intrinsic metric $\rho$) and the following condition on $V$:
\[
V(x)\geq -b_1-b_2[\rho(o,x)]^2,
\]
for all $x\in X$, where $o\in X$ is a fixed vertex, and $b_1\geq 0$ and $b_2\geq 0$ are constants.

As seen in corollary 2.3 of~\cite{Mil-18}, there is a similar phenomenon pertaining to Schr\"odinger oeprators on complete Riemannian manifolds. As an application of corollary~\ref{C:main-1}, example~\ref{E:ex-3} below describes a situation where our result guarantees the essential self-adjointness of the corresponding Schr\"odinger operator, while the criterion from proposition 2.2 in~\cite{Gol} is inconclusive. Lastly, we indicate that our theorem~\ref{T:main-1} also holds for magnetic Schr\"odinger operators. Nevertheless, in the interest of keeping the notations simpler, our formulation uses Schr\"odinger operators without magnetic fields.

The article consists of six sections, including the introduction. In section~\ref{S:main}, after describing the setting, notations, and operators, we state the main result (theorem~\ref{T:main-1}) and its corollary. In section~\ref{S:example} we provide an example illustrating corollary~\ref{C:main-1}. In section~\ref{S:prelim} we collect some auxiliary lemmas used in latter parts of the paper. The last two sections are devoted to the proofs of theorem~\ref{T:main-1} and corollary~\ref{C:main-1}.

\section{Main Result}\label{S:main}  In this section, we describe weighted graphs, the corresponding function spaces, and discrete Schr\"odinger operators. At the end of the section we state the main result.

\subsection{Weighted Graphs and Basic Function Spaces}\label{SS:setting}
By a \emph{symmetric weighted graph} $(X,b)$ we mean a countable \emph{vertex set} $X\neq \emptyset$ together with an \emph{edge weight function} $b\colon X\times X\to[0,\infty)$ satisfying the following conditions:
\begin{enumerate}
  \item [(b1)] $b (x, y) = b (y, x)$, for all $x,\,y\in X$;
  \item [(b2)] $b(x,x)=0$, for all $x\in X$;
  \item [(b3)] $\sum_{z\in X} b(x,z)<\infty$.
\end{enumerate}
We say that the graph $(X,b)$ is \emph{locally finite} if for all $x\in X$ we have $\displaystyle\sharp\,\{y\in X\colon b(x,y)>0\}<\infty$, where $\sharp\, G$ is understood as the number of elements in the set $G$.

We say that the two vertices $x,\, y\in X$  are \emph{connected by an edge} and write $x\sim y$ if $b(x, y) > 0$.  By a \emph{finite path} we mean a finite sequence of vertices $x_0,\,x_1,\,\dots,x_n$ such that $x_{j}\sim x_{j+1}$ for all $0\leq j\leq n-1$. We call the graph $(X,b)$ \emph{connected} if every two vertices $x,\,y\in X$ can be joined by a finite path.

We equip the set $X$ with the discrete topology. By $C(X)$ we denote the set of functions $f\colon X\to\mathbb{C}$, that is, continuous (complex-valued) functions with respect to the discrete topology. We use the symbol $C_{c}(X)$ to indicate the set of finitely supported elements of $C(X)$.

By a \emph{weight} we mean a function $\mu\colon X\to (0,\infty)$. Using the formula
\[
\mu(S):=\displaystyle\sum_{x\in S}\mu(x), \qquad S\subseteq X,
\]
the function $\mu$ gives rise to a Radon measure of full support. In this paper, we will use the two concepts interchangeably.

Using a weight $\mu$, we define the anti-duality pairing $(\cdot,\cdot)_{a}\colon C(X)\times C_{c}(X)\to \mathbb{C}$ by
\begin{equation}\label{E:a-d}
(u,w)_{a}:=\displaystyle\sum_{x\in X}\mu(x)u(x)\overline{w(x)}.
\end{equation}
This pairing induces an anti-linear isomorphism $w\mapsto (\cdot,w)_{a}$ between $C_{c}(X)$ and $C(X)'$. Here, $C(X)'$ is the continuous dual of $C(X)$, where the latter space is equipped with the topology of pointwise convergence.

The notation $\lw$ indicates the space of functions $f\in C(X)$ such that
\begin{equation}\label{E:l-p-def}
\displaystyle\sum_{x\in X}\mu(x)|f(x)|^2<\infty,
\end{equation}
where $|\cdot|$ is the modulus of a complex number.

The space $\lw$ is a Hilbert space with the inner product
\begin{equation}\label{E:inner-w}
(f,g):=\sum_{x\in X}\mu(x)f(x)\overline{g(x)}.
\end{equation}
The inner product $(\cdot,\cdot)$ in $\lw$ gives rise to the norm $\|\cdot\|$. Comparing with $(\cdot,\cdot)_{a}$, we see that if $f\in\lw$ and $g\in\xcomp$, then $(f,g)=(f,g)_{a}$.

By \emph{weighted vertex degree} $\textrm{Deg}(x)$  we mean
\begin{equation}\label{E:w-deg}
   \textrm{Deg}(x):=\frac{1}{\mu(x)}\sum_{y\in X}b(x,y), \quad x\in X.
\end{equation}
Lastly, as in section 11.4 of~\cite{KLW-book}, we define the squared gradient of a function $f\in C(X)$ :
\begin{equation}\label{E:grad-sq}
    |\nabla f|^2(x):=\frac{1}{\mu(x)}\sum_{y\in X}b(x,y)|f(x)-f(y)|^2,\quad x\in X.
\end{equation}

\subsection{Finiteness Condition (FC)}\label{SS:FC}
As in definition 3.16 in~\cite{MS-18}, we say that the triplet $(X,b,\mu)$ satisfies the \emph{finiteness condition}, abbreviated as (FC), if for all $x\in X$ the mapping $y\mapsto b(x,y)/\mu(y)$ belongs to $\lw$.

\subsection{Discrete Schr\"odinger Operators}\label{SS:expressions}
For a weighted graph $(X,b)$, a weight $\mu\colon X\to (0,\infty)$, and a function $V\colon X\to\mathbb{R}$ (called \emph{potential}), we now define the \emph{formal discrete Schr\"odinger operator} $\mathcal{L}_{V}\colon \mathcal{F}\to C(X)$. Here, the domain of $\mathcal{L}_{V}$ is described as
\begin{equation}\label{E:def-f}
\mathcal{F}:=\{f\in C(X)\colon \displaystyle\sum_{y\in X}b(x,y)|f(y)|<\infty\textrm{ for all } x\in X\},
\end{equation}
and the action of $\mathcal{L}_{V}$ is described as
\begin{equation}\label{E:magnetic-lap}
    \mathcal{L}_{V}f(x):=\frac{1}{\mu(x)}\sum_{y\in X}b(x,y)(f(x)-f(y))+V(x)f(x),\quad x\in X.
\end{equation}
(To make our notation simpler, we suppressed $\mu$ in the symbol $\mathcal{L}_{V}$.) The absolute convergence of the sum is guaranteed by the description of the domain $\mathcal{F}$. It turns out that $\mathcal{F}=C(X)$ if and only if $(X,b)$ is locally finite.

\subsection{Pseudo Metric} By \emph{pseudo metric} we mean a symmetric map $\rho\colon X\times X\to [0, \infty)$ with zero diagonal (that is, $\rho(x,x)=$ for all $x\in X$) and satisfying the (usual) triangle inequality.
\subsection{Condition (J)}\label{SS:jump}
By a \emph{jump size} $s$ of a pseudo metric $\rho$ we mean
\begin{equation}\label{E:jump-size}
s:=\sup\{\rho (x,y)\colon \textrm {for all }x,y\in X\textrm{ with }x\sim y\}.
\end{equation}
We say that a pseudo metric $\rho$ satisfies the condition (J) if the jump size $s$ of $\rho$ is finite.

\subsection{Intrinsic Pseudo Metric} We call a pseudo metric $\rho$ \emph{intrinsic metric} for a graph $(X,b,\mu)$ if
\begin{equation}\label{E:int-metric}
\sum_{y\in X}b(x,y)(\rho(x,y))^2\leq\mu(x),
\end{equation}
for all $x\in X$.

For a pseudo metric $\rho(\cdot,\cdot)$, a number $r\geq0$ and a point $o \in X$, the notation $B(o, r)$ indicates the set
\begin{equation}\label{E:b-pseudo}
B(o, r):=\{y\in X\colon \rho(o,y)\leq r\}.
\end{equation}

\subsection{Condition (B*)}\label{SS:B-star} Following section 11.4 of~\cite{KLW-book}, we say that $(X,b,\mu)$ with a pseudo metric $\rho$ satisfies the condition (B*) if the restriction of $\textrm{Deg}$ (weighted vertex degree as in~(\ref{E:w-deg})) to balls $B(o, r)$ (with respect to $\rho$ as in~(\ref{E:b-pseudo})) is bounded for all $o\in X$ and $r\geq 0$.

\begin{remark} As pointed out in section 11.4 of~\cite{KLW-book}, the condition (B*) follows from a (stronger) condition, which we label as (B) and which requires that the distance balls $B(o, r)$ be finite for all $o\in X$ and $r\geq 0$. On the other hand, as indicated in exercise 11.8 of~\cite{KLW-book}, the condition (B*) does not imply the condition (B). As shown in lemma 11.28 of~\cite{KLW-book}, the conditions (B) and (J) together force the graph to be locally finite.
\end{remark}

\subsection{Statement of the Main Result} Throughout this subsection we assume that $(X,b,\mu)$ satisfies the finiteness condition (FC).
By lemma 3.15 in~\cite{MS-18}, the property (FC) is equivalent to the inclusion $\mathcal{L}_{V}[\xcomp]\subseteq\lw$. This means that, in this context, we can consider $\mathcal{L}_{V}|_{\xcomp}$ as an operator in $\lw$. Moreover, in view of Green's formula (see lemma~\ref{L:L-1} below), for any real-valued function $V\in C(X)$, the operator ${L}_{V}|_{\xcomp}$ is symmetric (as an operator in $\lw$).

This brings us to the main result of our article.

\begin{thm}\label{T:main-1} Let $(X, b,\mu)$ be a weighted and connected graph satisfying the property (FC) as in section~\ref{SS:FC}. Furthermore, assume that there exists an intrinsic metric $\rho$ such that the condition (B*) is satisfied, where (B*) is as in section~\ref{SS:B-star}. Assume that $V=U-W$, where $U$ and  $W$ are functions such that $U(x)\geq 0$ and $W(x)\geq 0$ for all $x\in X$. Additionally, assume that there exist constants $c_1\geq 0$ and $c_2\geq 0$ such that
\begin{equation}\label{E:W-cond}
|\nabla W|^2(x) \leq c_1+c_2W(x),
\end{equation}
for all $x\in X$, where $|\nabla W|^2(x)$ is as in~(\ref{E:grad-sq}).

Then, $\mathcal{L}_{V}|_{\xcomp}$ is essentially self-adjoint.
\end{thm}

The following corollary of theorem~\ref{T:main-1} provides a more concrete essential self-adjointness result:

\begin{cor}\label{C:main-1} Let $(X, b,\mu)$ be a weighted and connected graph satisfying the property (FC) as in section~\ref{SS:FC}. Assume that there exists an intrinsic metric $\rho$ such that the conditions (B*) and (J) are satisfied, where (B*) and (J) are as in sections~\ref{SS:B-star} and~\ref{SS:jump} respectively. Furthermore, assume that $V\colon X\to\mathbb{R}$ satisfies the following property: there exists $o\in X$ and constants $b_1\geq 0$ and $b_2\geq 0$ such that
\begin{equation}\label{E:Q-cond}
V(x)\geq -b_1-b_2[\rho(o,x)]^2,
\end{equation}
for all $x\in X$.

Then, $\mathcal{L}_{V}|_{\xcomp}$ is essentially self-adjoint.
\end{cor}
\begin{rem} In example~\ref{E:ex-3} we describe a situation where corollary~\ref{C:main-1} guarantees the essential self-adjointness of $\mathcal{L}_{V}|_{\xcomp}$, while the criterion from proposition 2.2 in~\cite{Gol} is inconclusive.
\end{rem}

\section{An Example Pertaining to Corollary~\ref{C:main-1}}\label{S:example} Before describing our example, we recall the definition of the path (pseudo) metric on a graph. A pseudo metric $\rho=\rho_{\sigma}$ is called \emph{a path pseudo metric} if there exists a map $\sigma\colon X\times X\to [0,\infty)$ such that $\sigma(x,y)=\sigma(y,x)$, for all $x,\,y\in X$; $\sigma(x,y)>0$ if and only if $x\sim y$; and
\begin{align*}
\rho_{\sigma}(x,y)=\inf\{l_{\sigma}(\gamma)\colon &\gamma=(x_0,x_1,\dots,x_n),\, n\geq 1, \\
&\textrm{is a path connecting }x\textrm{ and }y\},
\end{align*}
where the length $l_{\sigma}$ of the path $\gamma=(x_0,x_1,\dots,x_n)$ is given by
\begin{equation}\label{E:l-sigma-def}
l_{\sigma}(\gamma)=\sum_{i=0}^{n-1}\sigma(x_i,x_{i+1}).
\end{equation}
Using this procedure, the function
\begin{equation}\label{E:sigma-def-klw}
\sigma(x,y):=\min\left\{(\textrm{Deg}(x))^{-1/2}, (\textrm{Deg}(y))^{-1/2}\right\},\qquad x\sim y,
\end{equation}
where $\textrm{Deg}(x)$ is as in~(\ref{E:w-deg}),
leads to the so-called~\emph{degree-path metric} $\rho_{\sigma}$. This terminology is consistent with definition 11.18 of~\cite{KLW-book}. According to lemma 11.19 in~\cite{KLW-book}, the degree-path metric $\rho_{\sigma}$ is intrinsic.

\begin{example}\label{E:ex-3} Consider a graph whose vertex set $X$  is arranged in a ``triangular" pattern so that the vertex $x_{1,1}$ is in the first row, vertices $x_{2,1}$ and $x_{2,2}$ are in the second row,  vertices $x_{3,1}$, $x_{3,2}$, and $x_{3,3}$ are in the third row, and so on.  The vertex $x_{1,1}$ is connected to $x_{2,1}$ and $x_{2,2}$. The vertex $x_{2,i}$, where $i=1,2$, is connected to every vertex $x_{3,j}$, where $j=1,2,3$. The pattern continues so that each of $k$ vertices in the $k$-th row is connected to each of $k+1$ vertices in the $(k+1)$-st row.

Let $\mu(x)=2k^{1/2}$ for every vertex $x$ in the $k$-th row, and let $b(x,y)\equiv 1$ for all vertices $x\sim y$.
Denoting
\[
\textrm{deg}(x):=\displaystyle\sum_{{y\in X}}b(x,y),
\]
we note that for all $k\geq 1$ and $j\geq 1$ we have $\textrm{deg}(x_{k,j})=2k$, and hence the weighted vertex degree from~(\ref{E:w-deg}) is
\begin{equation}\label{E:w-deg-calc}
\textrm{Deg}(x_{k,j})=\frac{\textrm{deg}(x_{k,j})}{\mu(x_{k,j})}=\frac{2k}{2k^{1/2}}=k^{1/2}.
\end{equation}
Thus, the formula~(\ref{E:sigma-def-klw}), applied to every vertex $x$ in the $k$-th row and every vertex $y$ in the $(k+1)$-st row, yields
\begin{align*}
&\sigma(x,y):=\min\left\{(\textrm{Deg}(x))^{-1/2}, (\textrm{Deg}(y))^{-1/2}\right\}\nonumber\\
&=\min\left\{k^{-1/4}, (k+1)^{-1/4}\right\}=(k+1)^{-1/4}.\nonumber
\end{align*}
Fixing $o=x_{1,1}$ and denoting by $\rho_{\sigma}$ the degree-path metric corresponding to $\sigma$ (as discussed in the prelude to this example), we have for all vertices $x$ in the $k$-th row, $k\geq 2$:
\[
\rho_{\sigma}(o,x)=\displaystyle\sum_{j=1}^{k-1}(j+1)^{-1/4}\geq \frac{k-1}{k^{1/4}}\geq \frac{k^{3/4}}{2}.
\]
Hence, for all vertices $x$ in the $k$-th row, $k\geq 2$:
\begin{equation}\label{E:ex-rho-est}
k^{3/2}\leq 4 [\rho_{\sigma}(o,x)]^2.
\end{equation}
For all vertices $x$ in the $k$-th row, $k\geq 1$, define $V(x):=-k^{1/2}$.

In view of~(\ref{E:ex-rho-est}), for all vertices $x$ in the $k$-th row, $k\geq 2$, we have
\[
V(x)=-k^{1/2}\geq -k^{3/2}\geq -4 [\rho_{\sigma}(o,x)]^2.
\]
For (the only) vertex $o=x_{1,1}$ in the first row, we have
\[
V(o)=-1=0-1=-4 [\rho_{\sigma}(o,o)]^2-1.
\]
Hence, the potential $V$ satisfies the condition~(\ref{E:Q-cond}) with $b_1=1$ and $b_2=4$.

We also note that $\rho_{\sigma}$ has a finite jump size. Additionally, being a degree-path metric, $\rho_{\sigma}$ is an intrinsic metric. Furthermore, the balls $B_{\rho_{\sigma}}(x,r)$, $x\in X$, $r\geq 0$, are finite; in particular, the condition (B*) is satisfied. Therefore, by corollary~\ref{C:main-1} the operator $\mathcal{L}_{V}$ is essentially self-adjoint on $\xcomp$.

Furthermore, it is easy to see that for every $c\in\mathbb{R}$, there exists a function $u\in\xcomp$ such that the inequality
\begin{equation}\nonumber
(\mathcal{L}_{V}u,u)\geq c\|u\|^2
\end{equation}
is not satisfied. Thus, the operator $\mathcal{L}_{V}|_{\xcomp}$ is not semi-bounded from below, which tells us that we cannot use corollary 11.5.5 from~\cite{MS-18}.

We now recall the self-adjointness result (specialized to operators without magnetic fields) of proposition 2.2 from~\cite{Gol}, which allows operators that are not lower semi-bounded:

\smallskip

\emph{Let $(X,b,\mu)$ be a weighted graph. Let $V\colon X\to\mathbb{R}$ and $\delta>0$. Take $\lambda\in\mathbb{R}$ so that
\begin{equation}\label{E:g-condition-1}
\{x\in X\colon \lambda+ \operatorname{Deg}(x)+V(x)=0\}=\emptyset,
\end{equation}
where $\operatorname{Deg}(x)$ is as in~(\ref{E:w-deg}) above.
Suppose that for every sequence of vertices $\{y_1,\,y_2,\dots\}$ such that $y_j\sim y_{j+1}$, $j\geq 1$, the following property holds:
\begin{equation}\label{E:g-condition-2}
\sum_{n=1}^{\infty}a_n\mu(y_n)=\infty,
\end{equation}
where
\[
a_n:=\prod_{j=1}^{n-1}\left(\left(\frac{\delta}{\operatorname{Deg}(y_j)}\right)^2+\left(1+\frac{\lambda+V(y_j)}{\operatorname{Deg}(y_j)}\right)^2\right), \quad n>1,
\]
and $a_1:=1$. Then $\mathcal{L}_{V}$ is essentially self-adjoint on $\xcomp$.}

We will show that proposition 2.2 from~\cite{Gol} is not applicable in our example.

Let $\lambda\in\mathbb{R}$ be such that~(\ref{E:g-condition-1}) is satisfied, with $V$ as in our example. Let $a_n$ be as in~(\ref{E:g-condition-2}) corresponding to the path $\gamma=(x_{1,1};\,x_{2,1};\, x_{3,1};\dots)$, the potential $V(x_{k,1})=-{k}^{1/2}$, $\delta>0$, and $\lambda$.

Then $a_1=1$, and upon remembering~(\ref{E:w-deg-calc}), for all $n\geq 2$ we have
\begin{equation}\nonumber
a_n\leq\prod_{k=1}^{n-1}\left(\frac{\delta}{k^{1/2}}+\left|1+\frac{\lambda-k^{1/2}}{k^{1/2}}\right|\right)^2=
\prod_{k=1}^{n-1}\frac{(\delta+|\lambda|)^2}{k}=\frac{(\delta+|\lambda|)^{2n-2}}{(n-1)!}.
\end{equation}

Therefore, keeping in mind that $\mu(x_{k,j})=2k^{1/2}$, we get
\[
\sum_{n=1}^{\infty}a_{n}\mu(x_{n,1})\leq 1+\sum_{n=2}^{\infty}\frac{2\sqrt{n}\cdot(\delta+|\lambda|)^{2n-2}}{(n-1)!}.
\]
By the ratio test, the series on the right hand side of this inequality converges. Hence, looking at~(\ref{E:g-condition-2}), we infer that proposition 2.2 from~\cite{Gol} cannot be used in our example.
\end{example}

\section{Preliminary Lemmas}\label{S:prelim}

We begin with an abstract fact from perturbation theory.

\subsection{Abstract Perturbation Result} The following result was established in theorem 5.7 of \cite{Okazawa-82}:

\begin{lem}\label{L:okazawa-5-7} Let $S$ and $G$ be nonnegative
symmetric operators in a Hilbert space $\mathscr{H}$ with inner
product $(\cdot,\cdot)_{\mathscr{H}}$ and norm $\|\cdot\|_{\mathscr{H}}$. Let $\mathscr{D}$ be
a linear subspace of $\mathscr{H}$ on which $S+G$ is essentially
self-adjoint (so that $\mathscr{D}$ is a core of the closure
$\overline{(S+G)|_{\mathscr{D}}}$ of $(S+G)|_{\mathscr{D}}$). Assume that the following
inequalities hold for all $u\in\mathscr{D}$:
\begin{equation}\label{E:assumption-S-1}
\|Su\|_{\mathscr{H}}+\|Gu\|_{\mathscr{H}}\leq a_1\|u\|_{\mathscr{H}}+a_2\|(S+G)u\|_{\mathscr{H}}
\end{equation}
and
\begin{equation}\label{E:assumption-S-2}
|\IM (Gu,Su)_{\mathscr{H}}|\leq \widetilde{a_1}\|u\|_{\mathscr{H}}^2+\widetilde{a_2}\|(S+G)u\|_{\mathscr{H}}\|u\|_{\mathscr{H}},
\end{equation}
where $a_1\geq 0$, $a_2\geq 0$, $\widetilde{a_1}\geq 0$ and $\widetilde{a_2}\geq 0$ are
constants. Then $S-G$ is essentially self-adjoint on $\mathscr{D}$.
\end{lem}

\subsection{Green's Formula} As in section 2.3.1 of~\cite{KLW-book}, for $f\in C(X)$ we define
\begin{equation}\label{E:difference}
\nabla_{x,y}f:=f(x)-f(y).
\end{equation}

In the sequel, we use the  following version of Green's formula; see proposition~ 1.5 in~\cite{KLW-book}:

\begin{lemma} \label{L:L-1}  Let $(X, b, \mu)$ be a weighted graph (not necessarily satisfying the property (FC)). Let $W\in C(X)$ be a real-valued function. Let $\mathcal{F}$ be as in~(\ref{E:def-f}), let $(\cdot,\cdot)_{a}$ be as in~(\ref{E:a-d}), and let $\nabla_{x,y}$ be as in~(\ref{E:difference}). Then, for all $f\in\mathcal{F}$ and $u\in\xcomp$, we have
\begin{align}\nonumber
&(\mathcal{L}_{W}f,u)_{a}=\sum_{x\in X}\mu(x)\mathcal{L}_{W}f(x)\overline{u(x)}=\sum_{x\in X}\mu(x)f(x)\overline{\mathcal{L}_{W}u(x)}\nonumber\\
&=\frac{1}{2}\displaystyle\sum_{x,y\in X}b(x,y)(\nabla_{x,y}f)\overline{\nabla_{x,y}u}+\displaystyle\sum_{x\in X}\mu(x)W(x)f(x)\overline{u(x)},\nonumber
\end{align}
with the sums converging absolutely.
\end{lemma}
\begin{remark}\label{R:R-green} If the graph $(X, b, \mu)$ satisfies the property (FC), then $\lw\subseteq\mathcal{F}$; see lemma 2.15 in~\cite{MS-18}.
\end{remark}
\subsection{Leibniz Rule} By lemma 2.25 in~\cite{KLW-book}, for  all $f,\, g\in C(X)$ we have
\begin{equation}\label{E:leib}
\nabla_{x,y}(fg)=f(x)(\nabla_{x,y}g)+(\nabla_{x,y}f)g(y)=g(x)(\nabla_{x,y}f)+(\nabla_{x,y}g)f(y).
\end{equation}

\section{Proof of Theorem~\ref{T:main-1}}\label{S:thm-1}

To simplify the notations, the maximal multiplication operator in $\lw$ corresponding to $F\in C(X)$ will be denoted by $F$, that is, $(Fu)(x):=F(x)u(x)$ with $\textrm{Dom}(F):=\{u\in \lw\colon Fu\in\lw\}$. In this proof $\mathcal{L}_{0}$ is as in~(\ref{E:magnetic-lap}) with $V\equiv 0$.

We will check that the conditions of lemma~\ref{L:okazawa-5-7} are satisfied for $S=\mathcal{L}_{U}$ and $G=W+1$.

Using the hypotheses $U\geq 0$ and $W\geq 0$ as well as lemma~\ref{L:L-1} we can see that $\mathcal{L}_{U}|_{\xcomp}$ and $(W+1)|_{\xcomp}$ are non-negative symmetric operators in $\lw$. Furthermore, keeping in mind the hypotheses on the graph $(X,b,\mu)$, the metric $\rho$, and the weighted vertex degree $\textrm{Deg}$, theorem 12.21 from~\cite{KLW-book} tells us that $\mathcal{L}_{U}+W+1$ is essentially self-adjoint on $\xcomp$.

We now show that the hypothesis~(\ref{E:assumption-S-1}) of lemma~\ref{L:okazawa-5-7} is satisfied for $S=\mathcal{L}_{U}|_{\xcomp}$ and $G=(W+1)|_{\xcomp}$.

With $c_1$ and $c_2$ as in~(\ref{E:W-cond}), we set $c_3:=c_1+c_2$. Then, for all $u\in\xcomp$ we have
\begin{align}
&\|(\mathcal{L}_{U}+W+1)u\|^2+c_3\|u\|^2\nonumber\\
&=\|\mathcal{L}_{U}u\|^2+\|(W+1)u\|^2+2\RE(\mathcal{L}_{U}u,(W+1)u)+c_3\|u\|^2\nonumber\\
&\geq \|\mathcal{L}_{U}u\|^2+\|(W+1)u\|^2+2\RE(\mathcal{L}_{0}u,(W+1)u)+c_3\|u\|^2,\nonumber
\end{align}
where the last inequality holds since $(Uu,(W+1)u)\geq 0$.

Thus, to verify that hypothesis~(\ref{E:assumption-S-1}) of lemma~\ref{L:okazawa-5-7} is fulfilled, it remains to show that
\begin{equation}\label{E:okazawa-1-ts}
2\RE(\mathcal{L}_{0}u,(W+1)u)\geq -c_3\|u\|^2,
\end{equation}
for all $u\in\xcomp$.

Writing
\[
W+1=\sqrt{W+1}\sqrt{W+1}
\]
and using lemma~\ref{L:L-1} along with~(\ref{E:leib}) we have
\begin{align}\label{E:1st}
&2\RE(\mathcal{L}_{0}u,(W+1)u)\nonumber\\
&=\RE \displaystyle\sum_{x,y\in X}b(x,y)(\nabla_{x,y}u)\nabla_{x,y}(\overline{u}\sqrt{W+1}\sqrt{W+1})\nonumber\\
&=\RE \displaystyle\sum_{x,y\in X}b(x,y)(\nabla_{x,y}u)\sqrt{W(x)+1}\nabla_{x,y}(\overline{u}\sqrt{W+1})\nonumber\\
&+\RE \displaystyle\sum_{x,y\in X}b(x,y)(\nabla_{x,y}u)(\nabla_{x,y}\sqrt{W+1})\overline{u(y)}\sqrt{W(y)+1},
\end{align}
for all $u\in\xcomp$.

For future use, we set
\begin{equation}\label{E:int-I}
I:=\RE \displaystyle\sum_{x,y\in X}b(x,y)(\nabla_{x,y}u)(\nabla_{x,y}\sqrt{W+1})\overline{u(y)}\sqrt{W(y)+1}.
\end{equation}

Keeping in mind~(\ref{E:leib}), we now rewrite the first term on the right hand side of the second equality in~(\ref{E:1st}) with the help of the formula
\[
(\nabla_{x,y}u)\sqrt{W(x)+1}=\nabla_{x,y}(u\sqrt{W+1})-(\nabla_{x,y}\sqrt{W+1})u(y),
\]
and return to~(\ref{E:1st}):
\begin{align}\label{E:2nd}
&2\RE(\mathcal{L}_{0}u,(W+1)u)\nonumber\\
&=\RE \displaystyle\sum_{x,y\in X}b(x,y)\nabla_{x,y}(u\sqrt{W+1})\nabla_{x,y}(\overline{u}\sqrt{W+1})\nonumber\\
&-\RE \displaystyle\sum_{x,y\in X}b(x,y)(\nabla_{x,y}\sqrt{W+1})u(y)\nabla_{x,y}(\overline{u}\sqrt{W+1})+I,
\end{align}
where $I$ is as in~(\ref{E:int-I}).

Keeping in mind~(\ref{E:leib}), we now rewrite the second term on the right hand side of~(\ref{E:2nd}) using the formula
\[
\nabla_{x,y}(\overline{u}\sqrt{W+1})=(\nabla_{x,y}\sqrt{W+1})\overline{u(x)}+(\nabla_{x,y}\overline{u})\sqrt{W(y)+1},
\]
and go back to~(\ref{E:2nd}):
\begin{align}\label{E:3rd}
&2\RE(\mathcal{L}_{0}u,(W+1)u)\nonumber\\
&=\displaystyle\sum_{x,y\in X}b(x,y)|\nabla_{x,y}(u\sqrt{W+1})|^2\nonumber\\
&-\RE \displaystyle\sum_{x,y\in X}b(x,y)(\nabla_{x,y}\sqrt{W+1})u(y)(\nabla_{x,y}\sqrt{W+1})\overline{u(x)}\nonumber\\
&-\RE \displaystyle\sum_{x,y\in X}b(x,y)(\nabla_{x,y}\sqrt{W+1})u(y)(\nabla_{x,y}\overline{u})\sqrt{W(y)+1}+I,
\end{align}
where $I$ is as in~(\ref{E:int-I}).

Looking~(\ref{E:int-I}), we see that the last two terms on the right hand side of~(\ref{E:3rd}) cancel each other out (because $\RE z=\RE\bar{z}$, for all $z\in\mathbb{C}$). Furthermore, we note that the first term on the right hand side of~(\ref{E:3rd}) is non-negative. Therefore, the estimate~(\ref{E:3rd}) leads to
\begin{align}\label{E:4th}
&2\RE(\mathcal{L}_{0}u,(W+1)u)\nonumber\\
&=\displaystyle\sum_{x,y\in X}b(x,y)|\nabla_{x,y}(u\sqrt{W+1})|^2\nonumber\\
&-\RE \displaystyle\sum_{x,y\in X}b(x,y)(\nabla_{x,y}\sqrt{W+1})u(y)(\nabla_{x,y}\sqrt{W+1})\overline{u(x)}\nonumber\\
&\geq -\RE \displaystyle\sum_{x,y\in X}b(x,y)(\nabla_{x,y}\sqrt{W+1})^2\overline{u(x)}u(y)\nonumber\\
&\geq -\frac{1}{2}\sum_{x,y\in X}b(x,y)(\nabla_{x,y}\sqrt{W+1})^2(|u(x)|^2+|u(y)|^2)\nonumber\\
&=-\displaystyle\sum_{x,y\in X}b(x,y)(\nabla_{x,y}\sqrt{W+1})^2|u(x)|^2,
\end{align}
for all $u\in\xcomp$, where in the second inequality we used the estimate
\[
\RE (\overline{u(x)}u(y))\leq |u(x)||u(y)|\leq \frac{|u(x)|^2}{2}+\frac{|u(y)|^2}{2},
\]
Note that the last equality in~(\ref{E:4th}) is true because the value of the sum with the coefficient $-\frac{1}{2}$ does not change if we interchange the roles of $x$ and $y$ in the corresponding summand.

Using the formula
\[
(\sqrt{\alpha}-\sqrt{\beta})^2=\frac{(\alpha-\beta)^2}{(\sqrt{\alpha}+\sqrt{\beta})^2},\qquad \alpha>0,\quad \beta>0,
\]
we have
\begin{align}\label{E:5th}
&(\nabla_{x,y}\sqrt{W+1})^2=(\sqrt{W(x)+1}-\sqrt{W(y)+1})^2\nonumber\\
&=\frac{(W(x)-W(y))^2}{(\sqrt{W(x)+1}+\sqrt{W(y)+1})^2}\nonumber\\
&\leq \frac{(W(x)-W(y))^2}{(\sqrt{W(x)+1})^2}=\frac{(\nabla_{x,y} W)^2}{W(x)+1}.
\end{align}

Combining~(\ref{E:5th}) with~(\ref{E:4th}) and remembering the definition~(\ref{E:grad-sq}), we have
\begin{align}\label{E:6th}
&2\RE(\mathcal{L}_{0}u,(W+1)u)\geq -\displaystyle\sum_{x,y\in X}b(x,y)(\nabla_{x,y}\sqrt{W+1})^2|u(x)|^2\nonumber\\
&\geq -\displaystyle\sum_{x\in X}\mu(x)|u(x)|^2(W(x)+1)^{-1}(\mu(x))^{-1}\displaystyle\sum_{y\in X} b(x,y)(\nabla_{x,y} W)^2\nonumber\\
&=-\displaystyle\sum_{x\in X}\mu(x)|u(x)|^2(W(x)+1)^{-1}|\nabla W|^2(x)\nonumber\\
&\geq -\displaystyle\sum_{x\in X}\mu(x)|u(x)|^2(W(x)+1)^{-1}(c_1+c_2W(x))\nonumber\\
&\geq -(c_1+c_2)\|u\|^2,
\end{align}
for all $u\in\xcomp$. Here, in the third  inequality we used the hypothesis~(\ref{E:W-cond}) and in the fourth inequality we used the estimate
\[
(W(x)+1)^{-1}(c_1+c_2W(x))\leq c_1+c_2.
\]
This proves~(\ref{E:okazawa-1-ts}) with $c_3=c_1+c_2$. Thus, the hypothesis~(\ref{E:assumption-S-1}) of lemma~\ref{L:okazawa-5-7} is satisfied for $S=\mathcal{L}_{U}|_{\xcomp}$ and $G=(W+1)|_{\xcomp}$.

We now check that the hypothesis~(\ref{E:assumption-S-2}) of lemma~\ref{L:okazawa-5-7} is satisfied for $S=\mathcal{L}_{U}|_{\xcomp}$ and $G=(W+1)|_{\xcomp}$.

Since $U$ and $W$ are real-valued functions, for all $u\in\xcomp$ we have
\begin{equation}\label{E:im-ref}
\IM ((W+1)u,\mathcal{L}_{U}u)=\IM (Wu,\mathcal{L}_{0}u),
\end{equation}
where $\mathcal{L}_{0}$ is as in~(\ref{E:magnetic-lap}) with $V=0$.

Using lemma~\ref{L:L-1} along with~(\ref{E:leib}) we have
\begin{align}\label{E:1st-im}
&2\IM(Wu,\mathcal{L}_{0}u)=\IM \displaystyle\sum_{x,y\in X}b(x,y)(\nabla_{x,y}(Wu))(\nabla_{x,y}\overline{u})\nonumber\\
&=\IM \displaystyle\sum_{x,y\in X}b(x,y)W(x)(\nabla_{x,y}u)(\nabla_{x,y}\overline{u})\nonumber\\
&+\IM \displaystyle\sum_{x,y\in X}b(x,y)u(y)(\nabla_{x,y}W)(\nabla_{x,y}\overline{u})\nonumber\\
&=\IM \displaystyle\sum_{x,y\in X}b(x,y)u(y)(\nabla_{x,y}W)(\nabla_{x,y}\overline{u}),
\end{align}
because the first term on the right hand side of the second equality is zero (as the imaginary part of a real number).

Combining~(\ref{E:1st-im}) with the estimate
\begin{align}\nonumber
&\left|b(x,y)u(y)(\nabla_{x,y}W)(\nabla_{x,y}\overline{u})\right|\nonumber\\
&\leq \frac{1}{2}\left(b(x,y)(\nabla_{x,y}W)^2|u(y)|^2+b(x,y)|\nabla_{x,y}u|^2\right),\nonumber
\end{align}
we obtain
\begin{align}\label{E:2nd-im}
&|\IM(Wu,\mathcal{L}_{0}u)|\nonumber\\
&\leq  \frac{1}{4}\displaystyle\sum_{x,y\in X}b(x,y)(\nabla_{x,y}W)^2|u(y)|^2+\frac{1}{4}\displaystyle\sum_{x,y\in X}b(x,y)|\nabla_{x,y}u|^2,\nonumber\\
&=\frac{1}{4}\displaystyle\sum_{x,y\in X}b(x,y)(\nabla_{x,y}W)^2|u(y)|^2+\frac{1}{2}(\mathcal{L}_{0}u,u),
\end{align}
for all $u\in\xcomp$, where in the equality we used lemma~\ref{L:L-1}.

Rewriting the first term on the right hand side of the equality in~(\ref{E:2nd-im}), remembering~(\ref{E:grad-sq}), and using the condition~(\ref{E:W-cond}), we have
\begin{align}\label{E:3rd-im}
&\frac{1}{4}\displaystyle\sum_{x,y\in X}b(x,y)(\nabla_{x,y}W)^2|u(y)|^2\nonumber\\
&=\frac{1}{4}\displaystyle\sum_{y\in X}|u(y)|^2\mu(y)(\mu(y))^{-1}\displaystyle\sum_{x\in X}b(x,y)(\nabla_{x,y}W)^2\nonumber\\
&=\frac{1}{4}\displaystyle\sum_{y\in X}|u(y)|^2\mu(y)|\nabla W|^2(y)\leq \frac{c_1}{4}\|u\|^2+\frac{c_2}{4}(Wu,u)\nonumber\\
&\leq \frac{c_1}{4}\|u\|^2+\frac{c_2}{4}((W+1)u,u),
\end{align}
for all $u\in\xcomp$.

Combining~(\ref{E:2nd-im}) and~(\ref{E:3rd-im}) we obtain
\begin{align}\label{E:4th-im}
&|\IM(Wu,\mathcal{L}_{0}u)|\nonumber\\
&\leq  \frac{c_1}{4}\|u\|^2+\frac{c_2}{4}((W+1)u,u)+\frac{1}{2}(\mathcal{L}_{0}u,u)\nonumber\\
&\leq  \frac{c_1}{4}\|u\|^2+\frac{c_2}{4}((W+1)u,u)+\frac{1}{2}(\mathcal{L}_{U}u,u)\nonumber\\
&\leq  \frac{c_1}{4}\|u\|^2+\frac{c_2+2}{4}((\mathcal{L}_{U}+W+1)u,u)\nonumber\\
&\leq  \frac{c_1}{4}\|u\|^2+\frac{c_2+2}{4}\|(\mathcal{L}_{U}+W+1)u\|\|u\|,
\end{align}
for all $u\in\xcomp$, where in the second inequality we used the assumption $U\geq 0$.

Remembering~(\ref{E:im-ref}), the estimate~(\ref{E:4th-im}) shows that the hypothesis~(\ref{E:assumption-S-2}) of lemma~\ref{L:okazawa-5-7} is satisfied for $S=\mathcal{L}_{U}|_{\xcomp}$ and $G=(W+1)|_{\xcomp}$. Hence, $\mathcal{L}_{U}-W-1$ is essentially self-adjoint on $\xcomp$, and the same is true for $\mathcal{L}_{U}-W=\mathcal{L}_{V}$. $\hfill\square$

\section{Proof of corollary~\ref{C:main-1}}
With $s$ as in~(\ref{E:jump-size}) and $V$ as in~(\ref{E:Q-cond}), define
\[
W(x):=b_1+b_2(\rho(o,x)+s)^2,\qquad U(x):=V(x)+b_1+b_2(\rho(o,x)+s)^2,\quad x\in X.
\]
Remembering that $s\geq 0$ and recalling~(\ref{E:Q-cond}) we get
\[
U(x)\geq V(x)+b_1+b_2[\rho(o,x)]^2\geq 0,
\]
for all $x\in X$.
We will now show that $W$ satisfies~(\ref{E:W-cond}).

We start with two preliminary observations. Using the triangle inequality for $\rho$ we have
\[
|\rho(o,x)-\rho(o,y)|\leq \rho(x,y),
\]
and
\[
\rho(o,y)\leq \rho(o,x)+ \rho(x,y)\leq \rho(o,x)+s,
\]
where the last inequality follows from~(\ref{E:jump-size}).

Looking at the definition of $W$ we have,
\begin{align}\label{E:cor-est-1}
&(W(x)-W(y))^2=[b_2(\rho(o,x)+s)^2-b_2(\rho(o,y)+s)^2]^2\nonumber\\
&=b_2^2(\rho(o,x)-\rho(o,y))^2(\rho(o,x)+\rho(o,y)+2s)^2\nonumber\\
&\leq b_2^2(\rho(x,y))^2(2\rho(0,x)+3s)^2\leq 9b_2^2(\rho(x,y))^2(\rho(0,x)+s)^2
\end{align}
where in the first inequality we used the two preliminary observations listed above.

We now refer to the definition~(\ref{E:grad-sq}) and estimate $|\nabla W|^2(x)$:
\begin{align*}
&|\nabla W(x)|^2=\frac{1}{\mu(x)}\sum_{y\in X}b(x,y)|W(x)-W(y)|^2\nonumber\\
&\leq\frac{1}{\mu(x)}\sum_{y\in X}b(x,y)(\rho(x,y))^29b_2^2(\rho(0,x)+s)^2\leq 9b_2^2(\rho(0,x)+s)^2,\nonumber
\end{align*}
where in the first inequality we used~(\ref{E:cor-est-1}) and in the second inequality we used the property~(\ref{E:int-metric}).

Thus, the condition~(\ref{E:W-cond}) is satisfied for $c_1=0$ and $c_2=9b_2$. Therefore, $U$ and $W$ satisfy the hypotheses of theorem~\ref{T:main-1}. From the definitions of $U$ and $W$ we see that $V=U-W$. Therefore, by theorem~\ref{T:main-1}, the operator $\mathcal{L}_{V}|_{\xcomp}$ is essentially self-adjoint. $\hfill\square$

\vskip 0.25in

\noindent\textbf{Data Availability} Data sharing is not applicable to this article as no data sets
were generated or analyzed during this study.

\vskip 0.25in

\noindent\textbf{Declarations}

\vskip 0.25in

\noindent\textbf{Conflict of Interest} The author has no conflict of interest to declare that
are relevant to the content of this article.

\end{document}